\documentclass{birkart}

\usepackage{amsmath}
\usepackage{amssymb}
\usepackage{amsfonts}

 \newtheorem{thm}{Theorem}[section]
 \newtheorem{cor}[thm]{Corollary}
 \newtheorem{lem}[thm]{Lemma}
 \newtheorem{prop}[thm]{Proposition}
 \theoremstyle{definition}

 \theoremstyle{remark}
 \newtheorem{rem}[thm]{Remark}

\numberwithin{equation}{section}
\numberwithin{figure}{section}

\newcommand{\hispace}{\mathcal H}
\newcommand{\gispace}{\mathcal G}
\newcommand{\nullspace}{\mathcal N}
\newcommand{\mospace}{\mathcal M}
\newcommand{\zeroset}{{\mathrm Z}}
\newcommand{\diff}{{\mathrm d}}

\newcommand{\Bergdiskgen}{A^2_{\alpha,\beta,\theta,\vartheta}(\mathbb{D}^2)}
\newcommand{\ZeroN}{\nullspace_{\alpha,\beta,\theta, \vartheta,N}(\D^2)}

\newcommand{\DiffN}
{\mospace_{\alpha,\beta,\theta,\vartheta,N}(\D^2)}

\newcommand{\Lball}{L^2_{\alpha,\beta,\theta}({\mathbb {B}}^2)}
\newcommand{\Aball}{A^2_{\alpha, \beta,\theta}({\mathbb B}^2)}
\newcommand{\Jball}{\mospace_{\alpha, \beta,\theta,N}({\mathbb B}^2)}
\newcommand{\LBargmann}{L^2_{\alpha,\beta,\theta}(\C^2)}
\newcommand{\Bargmann}{A^2_{\alpha,\beta,\theta}(\C^2)}
\newcommand{\ZeroNB}{\nullspace_{\alpha,\beta,\theta,N}(\C^2)}
\newcommand{\ZeroNpB}
{\nullspace_{\alpha,\beta,\theta,N+1}(\C^2)}
\newcommand{\DiffNB}
{\mospace_{\alpha,\beta,\theta,N}(\C^2)}

\newcommand{\C}{\mathbb C}
\newcommand{\D}{\mathbb D}

\newcommand{\Kerneldiskgen}
{P_{\alpha,\beta,\theta,\vartheta}}
\newcommand{\KerneldiskgenN}
{P_{\alpha,\beta,\theta+N,\vartheta}}
\newcommand{\KerneldiskN}
{P_{\alpha,\beta,\theta,\vartheta,N}}
\newcommand{\DiffkernelN}
{Q_{\alpha,\beta,\theta,\vartheta,N}}
\newcommand{\DiffkernelNpN}
{Q_{\alpha,\beta,\theta+N,\vartheta}}
\newcommand{\Diffkernelz}
{Q_{\alpha,\beta,\theta,\vartheta}}

\newcommand{\KernelB}{P_{\alpha,\beta,\theta}}
\newcommand{\KernelNB}
{P_{\alpha,\beta,\theta+N}}
\newcommand{\KernelBNB}
{P_{\alpha,\beta,\theta,N}}
\newcommand{\DiffkernelNB}
{Q_{\alpha,\beta,\theta,N}}
\newcommand{\DiffkernelNpNB}
{Q_{\alpha,\beta,\theta+N}}
\newcommand{\DiffkernelzB}
{Q_{\alpha,\beta,\theta}}

\newcommand{\diag}{\,{\rm diag}\,}

\begin{document}
\title[Norm expansion along a zero variety]
{Norm expansion along a zero variety in $\C^d$}
%----------Author 1
\author[Hedenmalm]
{H\aa{}kan Hedenmalm}

\address{Hedenmalm: Department of Mathematics\\
The Royal Institute of Technology\\
S -- 100 44 Stockholm\\
SWEDEN}

\email{haakanh@math.kth.se}

\thanks{Research supported by the G\"oran Gustafsson Foundation.}

%----------Author 2
\author[Shimorin]
{Serguei Shimorin}

\address{Shimorin: Department of Mathematics\\
The Royal Institute of Technology\\
S -- 100 44 Stockholm\\
SWEDEN}
\email{shimorin@math.kth.se}

%----------Author 3
\author[Sola]
{Alan Sola}

\address{Sola: Department of Mathematics\\
The Royal Institute of Technology\\
S -- 100 44 Stockholm\\
SWEDEN}
\email{alansola@math.kth.se}



%\date{March 23, 2004}
%\dedicatory{To my parents}
%%% ----------------------------------------------------------------------

\begin{abstract}
The reproducing kernel function of a weighted Bergman space over domains in 
$\C^d$ is known explicitly in only a small number of instances. 
Here, we introduce a process of orthogonal norm expansion along a subvariety 
of codimension $1$, which also leads to a series expansion of the reproducing 
kernel in terms of reproducing kernels defined on the subvariety. The 
problem of finding the reproducing kernel is thus reduced to the same kind of 
problem when one of the two entries is on the subvariety. A complete expansion
of the reproducing kernel may be achieved in this manner. We carry this out 
in dimension $d=2$ for certain classes of weighted Bergman spaces over the 
bidisk (with the diagonal $z_1=z_2$ as subvariety) and the ball 
(with $z_2=0$ as subvariety), as well as for a weighted Bargmann-Fock space
over $\C^2$ (with the diagonal $z_1=z_2$ as subvariety). 
\end{abstract}

%%% ----------------------------------------------------------------------
\maketitle
%%% ----------------------------------------------------------------------

%\addtolength{\oddsidemargin}{-1.3cm}
%\addtolength{\evensidemargin}{-1.3cm}
%\addtolength{\textwidth}{1.5cm}
%\addtolength{\textheight}{2.5cm}
%\addtolength{\topmargin}{-1.0cm}
%\addtolength{\footskip}{1.5cm}

%\begin{document}

%%% Research support

%\subjclass{}

%\keywords{}

%%% Section numbering

%\setcounter{section}{-1}
%\setcounter{equation}{0}
%\setcounter{thm}{0}
%\setcounter{prop}{0}
%\setcounter{lemma}{0}
%\setcounter{cor}{0}
%\setcounter{remark}{0}

\section{Introduction}

\noindent\bf The general setup. \rm 
Let $\Omega$ be an open connected set in $\C^d$ ($d=1,2,3,\ldots$). 
A separable Hilbert space $\hispace(\Omega)$ 
(over the complex field $\C$) of holomorphic functions on $\Omega$ is given, 
such that the point evaluations at points of $\Omega$ are bounded linear 
functionals on $\hispace(\Omega)$. By a standard result in Hilbert space 
theory, then, to each point $w\in\Omega$, there corresponds an element 
$k_w\in\hispace(\Omega)$ such that
$$f(w)=\langle f,k_w\rangle_{\hispace(\Omega)},\qquad f\in\hispace(\Omega).$$
Usually, we write $k(z,w)=k_w(z)$, and when we need to emphasize the space,
we write $k^{\hispace(\Omega)}$ in place of $k$. The function 
$k^{\hispace(\Omega)}$ is the {\em reproducing kernel} of $\hispace(\Omega)$. 
It is in general a difficult problem to calculate the reproducing kernel 
explicitly. Of course, in terms of an orthonormal basis $e_1,e_2,e_3,\ldots$ 
for $\hispace(\Omega)$, the answer is easy:
$$k(z,w)=\sum_{n=1}^{+\infty}e_n(z)\,\bar e_n(w).$$
In most situations where no obvious orthogonal basis is present, this 
requires application of the rather complicated Gram-Schmidt 
orthogonalization procedure. Here, we introduce a method which has the 
potential to supply the reproducing kernel in a more digestible form.
The method also supplies an expansion of the norm in $\hispace(\Omega)$ in
terms of norms of ``generalized restrictions'' along an analytic variety of 
codimension $1$. 

Let $p:\C^d\to\C$ be a nontrivial polynomial of $d$ variables, and let 
${\zeroset}_p$ be the variety 
$${\zeroset}_p=\big\{z\in\Omega:\,p(z)=0\big\},$$
which we assume to be nonempty. 
We also assume that $p$ has nonvanishing gradient along ${\zeroset}_p$.
This assures us that a holomorphic function in $\Omega$ that vanishes
on $\zeroset_p$ is analytically divisible by $p$ in $\Omega$. The assumptions 
made here are excessive, and may be relaxed substantially without 
substantially altering the assertions made in the sequel. For instance, 
$\Omega$ might instead be a $d$-dimensional complex manifold, and $p$ an
arbitrary analytic function on $\Omega$ with nonvanishing gradient along 
its zero set.
For $N=0,1,2,3,\ldots$, the subspace of $\hispace(\Omega)$ consisting 
functions holomorphically divisible in $\Omega$ by $p^N$ is denoted by
$\nullspace_N(\Omega)$; it is easy to show that $\nullspace_N(\Omega)$ is
a closed subspace of $\hispace(\Omega)$. We also need the difference space
$$\mospace_N(\Omega)=\nullspace_N(\Omega)\ominus\nullspace_{N+1}(\Omega),$$
which is a closed subspace of $\nullspace_N(\Omega)$. 
Let $P_N$ stand for the orthogonal 
projection $\hispace(\Omega)\to\nullspace_N(\Omega)$, while $Q_N$ is the
orthogonal projection $\hispace(\Omega)\to\mospace_N(\Omega)$. 
Let $\hispace_N(\Omega)$ be the Hilbert space of analytic functions $f$ on 
$\Omega$ such that $p^Nf\in\hispace(\Omega)$, with norm 
$$\|f\|_{\hispace_N(\Omega)}=\big\|p^Nf\big\|_{\hispace(\Omega)}.$$
Clearly, the operator $M_p^N$ of multiplication by $p^N$ is an isometric
isomorphism $\hispace_N(\Omega)\to\nullspace_N(\Omega)$. 
\medskip

\noindent\bf The norm expansion. \rm
We obtain a natural orthogonal decomposition
\begin{equation}
g=\sum_{N=0}^{+\infty}Q_N g,\qquad \|g\|^2_{\hispace(\Omega)}=
\sum_{N=0}^{+\infty}\|Q_Ng\|^2_{\hispace(\Omega)},
\label{norm-1}
\end{equation}
since
\begin{equation}
\bigcap_{N=0}^{+\infty}\nullspace_N(\Omega)=\{0\},
\label{emptyintersect}
\end{equation}
which expresses the fact that no analytic function on $\Omega$ may be 
holomorphically divisible by $p^N$ for all positive integers $N$ unless the 
function vanishes identically.  
In other words, we have an orthogonal decomposition
\[\hispace(\Omega)=\bigoplus_{N=0}^{+\infty}\mospace_{N}(\Omega).\] 
If we introduce the operator $R_N:\hispace(\Omega)
\to\hispace_N(\Omega)$ defined by $R_Ng=Q_Ng/p^N$, it is possible to write 
the above decomposition in the form 
$$g=\sum_{N=0}^{+\infty}p^N R_N g,\qquad \|g\|^2_{\hispace(\Omega)}=
\sum_{N=0}^{+\infty}\|R_Ng\|^2_{\hispace_N(\Omega)}.$$
The space of restrictions to $\zeroset_p$ of the functions in 
$\hispace_N(\Omega)$ is denoted by $\hispace_N(\zeroset_p)$. Is is supplied
with the induced Hilbert space norm
$$\|f\|_{\hispace_N({\zeroset}_p)}=\inf\big\{\|g\|_{\hispace_N(\Omega)}:\,g\in
\hispace_N(\Omega) \,\,\text{with}\,\, g|_{{\zeroset}_p}=f\big\}.$$
Let $\gispace_N(\Omega)$ denote the closed subspace of $\hispace_N(\Omega)$
consisting of $g$ with 
$p^Ng\in\mospace_N(\Omega)$. Also, let 
$\oslash_p$ denote the operation of taking the restriction to $\zeroset_p$
of a function defined on $\Omega$. It is easy to see that we have 
\begin{equation}
\|g\|_{\hispace_N(\Omega)}=
\big\|\oslash_p g\big\|_{\hispace_N({\zeroset}_p)}
\label{normid}
\end{equation}
if and only if $g\in\gispace_N(\Omega)$ (we recall that 
$g\in\hispace_N(\Omega)$ means that $p^Ng\in\nullspace_N(\Omega)$).
By polarizing (\ref{normid}), we find that 
\begin{equation}
\langle f,g\rangle_{\hispace_N(\Omega)}=
\big\langle\oslash_p f,\oslash_p g\big\rangle_{\hispace_N({\zeroset}_p)},
\qquad f,g\in\gispace_N(\Omega).
\label{scalid}
\end{equation}
Let $\oslash_p$ denote the operation of taking the restriction to $\zeroset_p$
of a function defined on $\Omega$.
We may now rewrite the orthogonal decomposition in yet another
guise (for $g\in\hispace(\Omega)$):
\begin{equation}
g=\sum_{N=0}^{+\infty}p^N R_N g,\qquad \|g\|^2_{\hispace(\Omega)}=
\sum_{N=0}^{+\infty}
\big\|\oslash_p R_Ng\big\|^2_{\hispace_N(\zeroset_p)}.
\label{normexp}
\end{equation}
In a practical situation, if we want to make use of this norm decomposition, 
we need to be able to characterize the restriction spaces 
$\hispace_N(\zeroset_p)$ in terms of a condition on $\zeroset_p$ (which has 
codimension $1$), and also to characterize the operators 
$\widetilde R_N=\oslash_p R_N$. This is quite often possible.
\medskip

\noindent\bf Expansion of the reproducing kernel. \rm
The above orthogonal decomposition corresponds to a reproducing kernel
decomposition
\begin{equation}
k^{\hispace(\Omega)}(z,w)=\sum_{N=0}^{+\infty}
k^{\mospace_N(\Omega)}(z,w)=
\sum_{N=0}^{+\infty}p(z)^N\bar p(w)^N k^{\gispace_{N}(\Omega)}(z,w).
\label{exp-1}
\end{equation}
Sometimes it is possible to characterize the restriction of 
$k^{\hispace_N(\Omega)}$ to $\Omega\times\zeroset_p$ (and hence, by symmetry,
to $\zeroset_p\times\Omega$ as well). One way this may happen is as follows.
Firstly, there is a certain point $w^0\in\zeroset_p$ for which it is easy to 
calculate the function $z\mapsto  k^{\hispace_{N}(\Omega)}(z,w^0)$ 
explicitly. Secondly, the automorphism group of $\Omega$ is fat enough, in
the sense that to each $w\in\zeroset_p$ there exists an automorphism of
$\Omega$ which sends $w^0$ to $w$. Moreover, to each automorphism we need
an associated unitary operator on $\hispace_{N}(\Omega)$ of composition type
(in more detail, it should be of the type $M_F C_\phi$, where 
$C_\phi f=f\circ \phi$ and $\phi$ is the automorphism in question, while 
$M_F$ denotes multiplication by a zero-free analytic function $F$).
The automorphisms allow us to calculate $k^{\hispace_{N}(\Omega)}(z,w)$
for $w\in\zeroset_p$ knowing $k^{\hispace_{N}(\Omega)}(z,w^0)$. Note that
on the set $\Omega\times\zeroset_p$, the two reproducing kernels 
$k^{\hispace_{N}(\Omega)}$ and $k^{\gispace_{N}(\Omega)}$ coincide.

Our goal is to express $k^{\hispace(\Omega)}$. 
Consider for a moment the following inner product:
\begin{equation}
l_N(z,w)=\big\langle \oslash_p k^{\hispace_{N}(\Omega)}_{w},
\oslash_p k^{\hispace_{N}(\Omega)}_{z}
\big\rangle_{\hispace_{N}(\zeroset_p)},\qquad (z,w)\in\Omega\times\Omega.
\label{def-l}
\end{equation}
Clearly, $l_N(z,w)$ is analytic in $z$ and antianalytic in $w$. 
Since $k^{\hispace_{N}(\Omega)}$ and $k^{\gispace_{N}(\Omega)}$ coincide on 
the set $\Omega\times\zeroset_p$, we have
$$l_N(z,w)=\big\langle \oslash_p k^{\gispace_{N}(\Omega)}_{w},
\oslash_p k^{\gispace_{N}(\Omega)}_{z}\big\rangle_{\hispace_{N}(\zeroset_p)},
$$
and if we apply (\ref{scalid}), we get 
\begin{equation}
l_N(z,w)=\big\langle k^{\gispace_{N}(\Omega)}_{w},
k^{\gispace_{N}(\Omega)}_{z}
\big\rangle_{\hispace_{N}(\Omega)}=k^{\gispace_{N}(\Omega)}(z,w).
\label{eq-l}
\end{equation}
By (\ref{exp-1}), we may now write down the desired explicit formula for 
$k^{\hispace(\Omega)}$, valid on $\Omega\times\Omega$:
\begin{equation}
k^{\hispace(\Omega)}(z,w)=
\sum_{N=0}^{+\infty}p(z)^N\bar p(w)^N 
\big\langle \oslash_p k^{\hispace_{N}(\Omega)}_{w},
\oslash_p k^{\hispace_{N}(\Omega)}_{z}\big\rangle_{\hispace_{N}(\zeroset_p)}.
\label{kernelexp}
\end{equation}

\medskip

\noindent\bf Applications. \rm
In Section \ref{Bidisk}, we carry out this program for classes of weighted 
Bergman spaces on the bidisk (with $p(z_1,z_2)=z_1-z_2$), while in Section
\ref{Ball}, we do the same thing for the ball in $\C^2$ 
(with $p(z_1,z_2)=z_2$). Finally, in Section \ref{Bargmann}, we apply the
technique to weighted Bargmann-Fock spaces on $\C^2$ 
(with $p(z_1,z_2)=z_1-z_2$). 

We remark that the first norm decomposition of this type for the bidisk 
was obtained by Hedenmalm and Shimorin \cite{HS}, who used it to
substantially improve the previously known estimates of the integral means 
spectrum for conformal maps.
\medskip

\noindent\bf A trivial example. \rm
Let $\diff A$ denote the normalized area element in the plane,
\begin{equation}
\diff A(z)=\frac{1}{\pi}\, 
\diff x\diff y, \hspace{0.2cm} \textrm{where}\hspace{0.2cm} z=x+iy,
\label{dA}
\end{equation}
and for $\alpha$, $-1<\alpha<+\infty$, we consider the following weighted area
element in the unit disk $\D=\{z\in\C:|z|<1\}$:
\begin{equation}
\diff A_\alpha(z)=(\alpha+1)(1-|z|^2)^\alpha\diff A(z).
\label{dAalpha}
\end{equation}
It is a probability measure in $\D$. The Hilbert space $A^2_\alpha(\D)$ 
consists of all analytic functions $g$ in $\D$ subject to the norm 
boundedness condition
\begin{equation}
\|g\|^2_{\alpha}=\int_\D|g(z)|^2\diff A_\alpha(z)<+\infty.
\label{A2alpha}
\end{equation}
Fix $\alpha$, and consider the space $\hispace(\D)=A^2_\alpha(\D)$ and
the polynomial $p(z)=z$. Then the space $\nullspace_N(\D)$ consists of
all functions that have a zero of order $N$ at the origin, while 
$\nullspace_N(\D)\ominus\nullspace_{N+1}(\D)$ is just the linear span of
the function $z^N$. We readily find that the orthogonal expansion
(\ref{normexp}) condenses to the familiar
$$g(z)=\sum_{N=0}^{+\infty}c_nz^n,\qquad \|g\|^2_\alpha=
\sum_{N=0}^{+\infty}\frac{N!}{(\alpha+2)_N}\,|c_N|^2,$$
where $(x)_n$ is the familiar Pochhammer symbol. The reproducing kernel for 
the space $A^2_\alpha(\D)$ is well-known:
$$k(z,w)=\sum_{N=0}^{+\infty}\frac{(\alpha+2)_N}{N!}\,z^N\bar w^N
=(1-z\bar w)^{-2-\alpha}.$$
The interesting thing is that the method outlined above applies to 
give the indicated representation of the reproducing kernel. 
\medskip

\noindent\bf Notation. \rm
In the rest of the paper, the notation for reproducing kernels is
slightly different (with letters $P$ and $Q$ instead of $k$). Also, we
should point out that in the sequel, the notation is consistent within each
section, but not necessarily between sections. This mainly applies to the 
spaces and their reproducing kernels, as we intentionally use very similar 
notation to demonstrate the analogy between the three cases we study (bidisk,
ball, Bargmann-Fock).  

\section{Weighted Bergman spaces in the bidisk}
\label{Bidisk}

\noindent\bf Preliminaries. \rm
The unit bidisk in $\mathbb{C}^2$ is the set
\[\D^2=\big\{z=(z_1,z_2) \in \mathbb{C}^2:\,|z_1|<1,\,|z_2|<1\big\}.\]
For a survey of the function theory of the bidisk, we refer to \cite{Rud1};
see also \cite{Krantz} and \cite{Berg}.
Fix real parameters $\alpha,\beta,\theta,\vartheta$ with 
$-1<\alpha,\beta,\theta,\vartheta<+\infty$.
We consider the Hilbert space 
$L^2_{\alpha,\beta,\theta,\vartheta}(\D^2)$ of
all (equivalence classes of) Borel measurable functions $f$ in the bidisk 
subject to the norm boundedness condition
\begin{equation*}
\|f\|^2_{\alpha,\beta,\theta,\vartheta}=
\int_{\D^2}|f(z_1,z_2)|^2|1-\bar{z_2}z_1|^{2\vartheta}
|z_1-z_2|^{2\theta}
\diff A_{\alpha}(z_1)\diff A_{\beta}(z_2)<+\infty,
%\label{lebspace}
\end{equation*}
where the notation is as in (\ref{dAalpha}); 
we let $\langle\cdot,\cdot\rangle_{\alpha,\beta,\theta,\vartheta}$ denote
the associated sesquilinear inner product. 
The {\em weighted Bergman space} $\Bergdiskgen$ is the subspace of 
$L^2_{\alpha,\beta,\theta,\vartheta}(\D^2)$
consisting of functions holomorphic in the bidisk. We need to impose a
further restriction on the parameters $\alpha,\beta,\theta,\vartheta$:
$$\alpha+\beta+2\theta+2\vartheta+3>0;$$ 
then the constant function $1$ will belong to the space $\Bergdiskgen$.

The reproducing kernel for $\Bergdiskgen$ will be denoted by
\[\Kerneldiskgen=
\Kerneldiskgen(z,w),\]
where we adhere to the notational convention
$$z=(z_1,z_2),\quad w=(w_1,w_2)$$
for points in $\C^2$. The kernel defines an orthogonal projection of 
the space $L^2_{\alpha,\beta,\theta,\vartheta}(\D^2)$ onto the weighted
Bergman space via the formula
\begin{multline}
P_{\alpha,\beta,\theta, \vartheta}[f](z)=\big\langle f,
\Kerneldiskgen
(\cdot,z)\big\rangle_{\alpha,\beta,\theta,\vartheta}\\
=\int_{\D^2}f(w)P_{\alpha,\beta,\theta,\vartheta}
(z,w)\,|1-\bar{w_2}w_1|^{2\vartheta}|w_1-w_2|^{2\theta}
\diff A_{\alpha}(w_1)\diff A_{\beta}(w_2);
\label{reprodprop}
\end{multline}
as indicated, we shall write $\Kerneldiskgen[f]$ for the projection of a 
function $f \in L^2_{\alpha,\beta,\theta,\vartheta}(\D^2)$.

In the case $\theta=\vartheta=0$, the reproducing kernel it is readily 
computed:
\begin{equation*}
P_{\alpha,\beta,0,0}(z,w)=
\frac{1}{(1-\bar{w_1}z_1)^{\alpha+2}(1-\bar{w_2}z_2)^{\beta+2}}.
\end{equation*}

We consider the polynomial $p(z_1,z_2)=z_1-z_2$ in the context of the 
introduction. In particular, for non-negative integers $N$, we consider 
the subspaces $\ZeroN$ of functions in $\Bergdiskgen$ that vanish up to 
order $N$ along the diagonal 
$$\diag(\D)=\big\{(z_1,z_2)\in \D^2:z_1=z_2\big\}.$$
Being closed subspaces of a reproducing kernel space, the spaces $\ZeroN$ 
possess reproducing kernels of their own. We shall write
\[\KerneldiskN=\KerneldiskN(z,w)\]
for these kernel functions.
The operators associated with the kernels project the space
$L^2_{\alpha,\beta,\theta,\vartheta}(\D^2)$
orthogonally onto $\ZeroN$. As before, we write
$\KerneldiskN[f]$ for the projection of a function.

Next, we define the spaces $\DiffN$ by setting
\[\DiffN=\ZeroN \ominus \ZeroNpB.\]
The spaces $\DiffN$ also admit reproducing kernels, and their
kernel functions are of the form
\begin{equation*}\DiffkernelN(z,w)
=\KerneldiskN(z,w)
-P_{\alpha,\beta,\theta,\vartheta,N+1}(z,w).
\end{equation*}
We shall write $\Diffkernelz$ for the kernel 
$Q_{\alpha,\beta,\theta,\vartheta,0}$.
\medskip

We begin with the following observation.

\begin{lem}
We have
\begin{equation}\KerneldiskN(z,w)
=(z_1-z_2)^N(\bar{w_1}-\bar{w_2})^N
\KerneldiskgenN(z,w)
\label{kerrel1}
\end{equation}
for $z,w\in\C^2$. 
\label{lemma-1}
\end{lem}

\begin{proof}
After multiplying both sides of (\ref{reprodprop}) by $(z_1-z_2)^N$ and using 
the fact that $|w_1-w_2|^{2N}=(w_1-w_2)^{N}(\bar{w_1}-\bar{w_2})^{N}$, 
we see that 
\begin{multline*}(z_1-z_2)^Nf(z)=
\int_{\D^2} \big\{(z_1-z_2)^N(\bar{w_1}-\bar{w_2})^N
\KerneldiskgenN(z,w)\big\}\\
\times \big\{(w_1-w_2)^Nf(w)\big\}\,|1-\bar{w_2}w_1|^{2\vartheta}
|w_1-w_2|^{2\theta}\diff A_{\alpha}(w_1)\diff A_{\beta}(w_2)
\end{multline*}
for every $f \in A^2_{\alpha,\beta,\theta+N,\vartheta}(\D^2)$.
From this it follows that 
$(z_1-z_2)^N(\bar{w_1}-\bar{w_2})^N\KerneldiskgenN(z,w)$
has the reproducing property for the space $\ZeroN$, and the proof is 
complete.
\end{proof}

If we write, as in the introduction, 
$$\hispace(\D^2)=A^2_{\alpha,\beta,\theta,\vartheta}(\D^2),$$
the argument of the proof of Lemma \ref{lemma-1} actually shows that we
have identified the spaces $\hispace_N(\D^2)$,
$$\hispace_N(\D^2)=A^2_{\alpha,\beta,\theta+N,\vartheta}(\D^2),
\qquad N=0,1,2,\ldots.$$
At the same time, we have also identified the spaces $\gispace_N(\D^2)$,
$$\gispace_N(\D^2)=
\mospace_{\alpha,\beta,\theta+N,\vartheta,0}(\D^2),\qquad N=0,1,2,\ldots.$$
As a consequence, we get that
\begin{equation}
\DiffkernelN(z,w)
=(z_1-z_2)^N(\bar{w_1}-\bar{w_2})^N
\DiffkernelNpN(z,w).
\label{kerrel2}
\end{equation}
By (\ref{exp-1}), we have the kernel function expansion
\begin{equation}
\Kerneldiskgen(z,w)
=\sum_{N=0}^{+\infty}\DiffkernelN(z,w),
\label{kerneldecomp}
\end{equation}
while the orthogonal norm expansion (\ref{norm-1}) reads
\begin{equation}
\|f\|^2_{\alpha,\beta,\theta,\vartheta}
=\sum_{N=0}^{+\infty}\big\|\DiffkernelN[f]
\big\|^2_{\alpha,\beta,\theta,\vartheta},
\qquad f\in\Bergdiskgen. 
\label{normdecomp}
\end{equation}

Our next objective is to identify the Hilbert space of restrictions to the
diagonal of 
$\hispace_N(\D^2)=A^2_{\alpha,\beta,\theta+N,\vartheta}(\D^2)$, as
well as to calculate the reproducing kernel of $\hispace_N(\D^2)$ on the
set $\D^2\times\diag(\D)$.
\medskip

\noindent\bf Unitary operators.
\rm 
The rotation operator $R_{\phi}$ (for a real parameter $\phi$) defined for 
$f\in\Bergdiskgen$ by
\[R_{\phi}[f](z_1,z_2)=f(e^{i\phi}z_1,e^{i\phi}z_2)\]
is clearly unitary, and we shall make use of it shortly. 
The following lemma supplies us with yet another family of unitary operators.

\begin{lem}
For every $\lambda \in \D$, the operator
\begin{equation}
U_{\lambda}[f](z_1,z_2)=
\frac{(1-|\lambda|^2)^{\alpha/2+\beta/2+\theta+\vartheta+2}}
{(1-\bar{\lambda}z_1)^{\alpha+\theta+\vartheta+2}
(1-\bar{\lambda}z_2)^{\beta+\theta+\vartheta+2}}
f\left(\frac{\lambda-z_1}{1-\bar{\lambda}z_1},
\frac{\lambda-z_2}{1-\bar{\lambda}z_2}\right)
\label{unitop}
\end{equation}
is unitary on the space $\Bergdiskgen$, and 
$U^2_{\lambda}[f]=f$ holds for every $f\in \Bergdiskgen$.
\end{lem}

\begin{proof}
For real parameters $p$ and $q$, we define the operator
\[\widetilde{U}_{\lambda}[f](z_1,z_2)
=\left(\frac{1-|\lambda|^2}{(1-\bar{\lambda}z_1)^2}\right)^p
\left(\frac{1-|\lambda|^2}{(1-\bar{\lambda}z_2)^2}\right)^q
f\left(\frac{\lambda-z_1}{1-\bar{\lambda}z_1},\frac{\lambda-z_2}
{1-\bar{\lambda}z_2}\right).\]
We want to choose $p$ and $q$ so that $\widetilde{U}_{\lambda}$ becomes
unitary.
A change of variables shows that
\begin{multline*}
\int_{\D^2}\left|\widetilde{U}_{\lambda}[f](z_1,z_2)\right|^2
|1-\bar{z_2}z_1|^{2\vartheta}|z_1-z_2|^{2\theta}
\diff A_{\alpha}(z_1)\diff A_{\beta}(z_2)\\
=\int_{\D^2}|f(\zeta,\xi)|^2
\frac{(1-|\lambda|^2)^{\alpha+\beta+2\theta+2\vartheta+4-2(p+q)}}
{|1-\bar{\lambda}\zeta|^{2\alpha+2\theta+2\vartheta+4-4p}
|1-\bar{\lambda}\xi|^{2\beta+2\theta+2\vartheta+4-4q}}
\,|1-\bar{\xi}\zeta|^{2\vartheta}
|\zeta-\xi|^{2\theta}\diff A_{\alpha}(\zeta)\diff A_{\beta}(\xi)
\end{multline*}
and we see that $p=1+(\alpha+\theta+\vartheta)/2$ and
$q=1+(\beta+\theta+\vartheta)/2$ are the correct choices. 

The proof of the second assertion is straight\-forward and therefore 
omitted.   
\end{proof}
\medskip

\noindent\bf The reproducing kernel on the diagonal. \rm
We now use the operators $R_{\phi}$ and $U_{\lambda}$ to 
compute the reproducing kernel on the set $\D^2\times\diag(\D)$.
We recall the standard definition of the generalized Gauss hypergeometric 
function
$${}_3 F_2\bigg(\begin{array}{ccc}a_1&a_2&a_3\\{}&b_1&b_2
\end{array}\bigg|x\bigg)=1+\sum_{n=1}^{+\infty}
\frac{(a_1)_n(a_2)_n(a_3)_n}{(b_1)_n(b_2)_n\,n!}\,x^n.$$

\begin{thm} We have that
\begin{equation*}
\Kerneldiskgen((z_1,z_2),(w_1,w_1))=\Diffkernelz((z_1,z_2),(w_1,w_1))=
\frac{\sigma(\alpha,\beta,\theta,\vartheta)}
{(1-\bar{w_1}z_1)^{\alpha+\theta+\vartheta+2}
(1-\bar{w_1}z_2)^{\beta+\theta+\vartheta+2}}
%\label{diagker}
\end{equation*}
for $z_1,z_2,w_1\in\D$. Here, $\sigma(\alpha,\beta,\theta,\vartheta)$ is the 
positive constant  given by 
\begin{multline*}
\frac{1}{\sigma(\alpha,\beta,\theta,\vartheta)}
=\int_\D\int_{\D}|1-\bar{z_2}z_1|^{2\vartheta}|z_1-z_2|^{2\theta}
\diff A_{\alpha}(z_1)\diff A_{\beta}(z_2).\\
=\frac{(\beta+1)\Gamma(\alpha+2)\Gamma(\theta+1)}
{(\alpha+\beta+2\theta+2\vartheta+3)\Gamma(\alpha+\theta+2)}\,\,
{}_3F_2\left(\begin{array}{ccc}
\theta+1 & \alpha+\theta+\vartheta+2 & \alpha+\theta+\vartheta+2 
\\ {}&\alpha+\theta+2 & \alpha+\beta+2\theta+2\vartheta+4 
\end{array}\bigg|1\right). 
%\label{sigma}
\end{multline*}
\label{thm-1}
\end{thm}

%!!!!alpha och beta boer kunna bytas ut mot varann utan att sigma aendras. Sant????

\begin{proof}
By the reproducing property of $\Kerneldiskgen$, we have 
\[f(0)=\big\langle f,\Kerneldiskgen(\cdot,0)
\big\rangle_{\alpha,\beta,\theta,\vartheta},\]
where $0$ this time denotes the origin in $\C^2$. The unitarity of 
$R_{\phi}$ gives us 
\[f(0)=\big\langle 
R_{\phi}[f],R_{\phi}[\Kerneldiskgen]
(\cdot,0)\big\rangle_{\alpha,\beta,\theta,\vartheta},\]
and since $f(0,0)=R_{\phi}[f](0,0)$, we see from the uniqueness of the 
reproducing kernel that
\[\Kerneldiskgen((e^{i\phi}z_1,e^{i\phi}z_2),(0,0))
=\Kerneldiskgen((z_1,z_2),(0,0))=\Kerneldiskgen(z,0).\]
The function $\Kerneldiskgen(z,0)$ is holomorphic in 
$\D^2$ and can be expanded in a power series. After comparing the 
series expansion  for the expressions on both sides of the above equality, 
we conclude that $\Kerneldiskgen(z,0)$ must be a (positive) constant,
which we denote by $\sigma(\alpha,\beta,\theta,\vartheta)$. 

Next, take $\lambda \in \D$ and $f \in \Bergdiskgen$.
Since the operators $U_{\lambda}$ are unitary and since $U^2_{\lambda}[f]=f$ 
we obtain that
\begin{multline*}
(1-|\lambda|^2)^{\alpha/2+\beta/2+\theta+\vartheta+2}
f(\lambda,\lambda)
=U_{\lambda}[f](0)
=\big\langle U_{\lambda}[f],
\Kerneldiskgen(\cdot,0)
\big\rangle_{\alpha,\beta,\theta,\vartheta}\\
=\big\langle U_{\lambda}^2[f], U_{\lambda}[\Kerneldiskgen
(\cdot,0)]\big\rangle_{\alpha,\beta,\theta,\vartheta}
=\sigma(\alpha,\beta,\theta)\,
\big\langle f,U_{\lambda}[1]\big\rangle_{\alpha,\beta,\theta,\vartheta}.
\end{multline*}
This equality together with the uniqueness of reproducing kernels
establishes that
\[\Kerneldiskgen((z_1,z_2),(\lambda,\lambda))
=\sigma(\alpha,\beta,\theta,\vartheta)
(1-|\lambda|^2)^{-\alpha/2-\beta/2-\theta-\vartheta-2}\,
U_{\lambda}[1](z_1,z_2),\]
which is the desired result.
The explicit expression for the constant in terms of an integral over the 
bidisk follows if we apply the reproducing property of the kernel applied 
to the constant function $1$; the evaluation of the integral in terms
of the hypergeometric function is done by performing the change of
variables 
\[z_1=\frac{z_2-\zeta}{1-\bar{z_2}\zeta},\quad z_2=z_2,\]
and by carrying out some tedious but straightforward calculations.
\end{proof}
\medskip

\noindent\bf Restrictions of reproducing kernels. \rm
From the previous subsection, we have that 
\[\Kerneldiskgen((z_1,z_2),(w_1,w_1))
=\frac{\sigma(\alpha,\beta,\theta,\vartheta)}
{(1-\bar{w_1}z_1)^{\alpha+\theta+\vartheta+2}
(1-\bar{w_1}z_2)^{\beta+\theta+\vartheta+2}}.\]
For continuous functions 
$f \in L^2_{\alpha,\beta,\theta,\vartheta}(\D^2)$, 
we use the notation $\oslash f$ for the restriction 
to the diagonal of the function, that is,
\[(\oslash f) (z_1)=f(z_1,z_1).\]
We fix $(w_1,w_1)$ and 
apply this operation to the kernel function of $\Bergdiskgen$. 
We obtain
\[\oslash\, \Kerneldiskgen(z_1,(w_1,w_1))
=\frac{\sigma(\alpha,\beta,\theta,\vartheta)}
{(1-\bar{w_1}z_1)^{\alpha+\beta+2\theta+2\vartheta+4}}\]
and we see that the restriction of the kernel
coincides with a multiple of the kernel function for
the space $A^2_{\alpha+\beta+2\theta+2\vartheta+2}(\D)$.
By the theory of reproducing kernels (see \cite{Sa}), this means that 
the induced norm for the space 
$\mospace_{\alpha,\beta,\theta,\vartheta,0}(\D^2)$ coincides with a 
multiple of the norm in the aforementioned 
weighted Bergman space in the unit disk. 
An immediate consequence of this fact is the inequality
\begin{equation}
\frac{1}{\sigma(\alpha,\beta,\theta, \vartheta)}
\|\oslash f\|^2_{\alpha+\beta+2\theta+2\vartheta+2}
\leq \|f\|^2_{\alpha,\beta,\theta,\vartheta},\qquad f\in\Bergdiskgen,
\label{rnorineq}
\end{equation}
and, more importantly, the equality
\begin{equation}
\frac{1}{\sigma(\alpha,\beta,\theta, \vartheta)}
\|\oslash f\|^2_{\alpha+\beta+2\theta+2\vartheta+2}
=\big\|\Diffkernelz[f]\big\|^2_{\alpha,\beta,\theta,\vartheta},
\qquad f\in\mospace_{\alpha,\beta,\theta,\vartheta,0}(\D^2).
\label{rnormeq}
\end{equation}
The notation on the left hand sides of (\ref{rnorineq}) and (\ref{rnormeq})
is in conformity with (\ref{A2alpha}).

The next step in our program is to compute the kernel function 
$\Diffkernelz$. 
In fact, we can determine the kernel in terms of an integral 
formula.

\begin{lem}
The kernel function for the space 
$\mospace_{\alpha,\beta,\theta,\vartheta,0}(\D^2)$ is given by 
\begin{equation*}
\Diffkernelz(z,w)
=\int_{\D}
\frac{\sigma(\alpha,\beta,\theta,\vartheta)\,\diff A_{\alpha+\beta+2\theta
+2\vartheta+2}(\xi)}
{[(1-\bar{\xi}z_1)(1-\xi \bar{w_1})]^{\alpha+\theta+\vartheta+2}
[(1-\bar{\xi}z_2)(1-\xi \bar{w_2})]^{\beta+\theta+\vartheta+2}}
%\label{intker}
\end{equation*}
for $z,w\in \D^2$.
\label{lemma-2}
\end{lem}

\begin{proof}
In the notation of the introduction, this is the identity (with $N=0$)
$$k^{\gispace_{N}(\Omega)}(z,w)=
\big\langle \oslash_p k^{\hispace_{N}(\Omega)}_{w},
\oslash_p k^{\hispace_{N}(\Omega)}_{z}
\big\rangle_{\hispace_{N}(\zeroset_p)},$$
which follows from (\ref{def-l}) and (\ref{eq-l}).
\end{proof}

We may now replace the terms on the right hand side in (\ref{normdecomp}) 
by norms taken in weighted spaces in the unit disk. 

\begin{lem}
For each $N=0,1,2,\ldots$, we have the equality of norms 
\begin{equation*}
\big\|\DiffkernelN[f]\big\|^2_{\alpha,\beta,\theta,\vartheta}
=\frac{1}{\sigma(\alpha,\beta,\theta+N,\vartheta)}
\left|\left|\oslash \left[\frac{P_{\alpha,\beta,\theta,\vartheta,N}[f]}
{(z_1-z_2)^N}\right]\right|\right|^2_{\alpha+\beta+2\theta+2\vartheta+2N+2},
\end{equation*}
for all $f\in\Bergdiskgen$.
\end{lem}

\begin{proof}
The statement follows from a combination of (\ref{kerrel1}) and 
(\ref{rnormeq}).
\end{proof}

We need one more result in order to complete the norm expansion for the 
bidisk.
In what follows, we use the notation $\partial_{z_1}f$ for the partial 
derivative of $f$ with respect to the variable $z_1$.

\begin{lem}
For $k=0,1,2,\ldots$, we have, for each $f\in\Bergdiskgen$,
\begin{equation*}
\oslash\, [\partial_{z_1}^{k}f]=
\sum_{n=0}^{k}n!{k\choose n}
\frac{(\alpha+\theta+\vartheta +n+2)_{k-n}}{
(\alpha+\beta+2\theta+2\vartheta+2n+4)_{k-n}}\,\partial_{z_1}^{k-n}
\oslash \left[\frac{P_{\alpha,\beta,\theta,\vartheta,n}[f]}{(z_1-z_2)^n}
\right].
\end{equation*}
\label{lemma-3}
\end{lem}

\begin{proof}
We recall that
\[\Kerneldiskgen (z,w)=
\sum_{N=0}^{+\infty}(z_1-z_2)^N(\bar{w_1}-\bar{w_2})^N
Q_{\alpha,\beta,\theta+N,\vartheta}(z,w),\]
whence it follows that
\begin{equation}
\oslash\, \partial_{z_1}^k\Kerneldiskgen(z_1,(w_1,w_2))=
\sum_{n=0}^k n!{k\choose n}(\bar{w_1}-\bar{w_2})^N
\oslash\, \partial_{z_1}^{k-n}Q_{\alpha,\beta,\theta+n,\vartheta}
(z_1,(w_1,w_2)).
\label{slashdiffrel}
\end{equation}
Differentiation of the integral formula of Lemma \ref{lemma-2} and taking the 
diagonal restriction leads to the equality
\begin{multline*}
\oslash\, \partial_{z_1}^{k-n}Q_{\alpha,\beta,\theta+n,\vartheta}
(z_1,(w_1,w_2))
=(\alpha+\theta+\vartheta+n+2)_{k-n}\,
\sigma(\alpha,\beta,\theta+n,\vartheta)\\
\times \int_{\D}\frac{\bar{\xi}^{k-n}\diff A_{\alpha+\beta+2\theta+
2\vartheta+2n+2}(\xi)}{(1-\bar{\xi}z_1)^{\alpha+\beta+2\theta+2\vartheta
+2n+4-(k-n)}(1-\xi \bar{w_1})^{\alpha+\theta+\vartheta+n+2}
(1-\xi \bar{w_2})^{\beta+\theta+\vartheta+n+2}}.
\end{multline*}
We now note that the expression
\[\frac{\bar\xi^{k-n}}
{(1-\bar{\xi}z_1)^{\alpha+\beta+2\theta+2\vartheta+2n+4+(k-n)}}\]
is a multiple of the reproducing kernel of the space
$A^2_{\alpha+\beta+2\theta+2\vartheta+2}(\D)$, differentiated 
$k-n$ times.
Invoking the reproducing property of this kernel, we obtain 
that
\[\oslash\, \partial_{z_1}^{k-n}Q_{\alpha,\beta,\theta+n,\vartheta}
(z_1,(w_1,w_2))=\frac{(\alpha+\theta+\vartheta+n+2)_{k-n}}
{(\alpha+\beta+2\theta+2\vartheta+2n+4)_{k-n}}\,\partial_{z_1}^{k-n}\oslash
\big[P_{\alpha,\beta,\theta+n,\vartheta}\big](z_1,(w_1,w_2)).\]
This result, together with the identities (\ref{slashdiffrel}) and
(\ref{kerrel1}), yields the desired equality, and the proof is complete.
\end{proof}
We remark that Lemma \ref{lemma-3} is rather the opposite to what we need;
it expresses the known quantity $\oslash[\partial_{z_1}^k f]$ in terms of
the quantities we should like to understand. Nevertheless, it is possible to
invert the assertion of Lemma \ref{lemma-3} and express the unknown quantities
in terms of known quantities. 
\medskip

\noindent\bf The diagonal norm expansion for the bidisk. \rm 
The above lemma finally allows us to express each term in the right-hand side 
of (\ref{normdecomp}) in terms of diagonal restrictions of derivatives of
the original function.

\begin{lem}
Put
\begin{equation*}
a_{k,N}=\frac{(-1)^{N-k}}{k!(N-k)!}\frac{(\alpha+\theta+\vartheta+k+2)_{N-k}}
{(\alpha+\beta+2\theta+2\vartheta+N+k+3)_{N-k}}.
\end{equation*}
Then, for all $N=0,1,2,\ldots$, the equality
\begin{equation*}
\oslash \left[ \frac{P_{\alpha,\beta,\theta,\vartheta}[f]}{(z_1-z_2)^N}
\right] 
=\sum_{k=0}^{N}a_{k,N}\,\partial_{z_1}^{N-k}\oslash \left[
\partial_{z_1}^k f\right],
\end{equation*}
holds for each $f \in \Bergdiskgen$.
\end{lem}

\begin{proof}
In view of the previous lemma it is enough to check that
\[\sum_{k=n}^Na_{k,N}n!{k\choose n}
\frac{(\alpha+\theta+\vartheta+n+2)_{k-n}}{(\alpha+\beta+2\theta
+2\vartheta+2n+4)_{k-n}}=\delta_{n,N},\]
where $\delta_{n,N}$ is the Kronecker delta. This 
amounts to performing some rather 
straight-forward calculations.
First we note that
\[(\alpha+\theta+\vartheta+k+2)_{N-k}(\alpha+\theta+\vartheta+n+2)_{k-n}
=(\alpha+\theta+\vartheta+2+n)_{N-n}\]
and since this last expression does not depend on $k$, we can factor it 
out from the above sum. This reduces our task to showing that 
\begin{equation*}
\sum_{k=n}^N\frac{(-1)^{N-k}}{(N-k)!(k-n)!
(\alpha+\beta+2\theta+2\vartheta+N+k+3)_{N-n}
(\alpha+\beta+2\theta+2\vartheta+2n+4)_{k-n}}=\delta_{n,N}.
\end{equation*}
We note that this is true when $n=N$. It remains to show that the left hand
side vanishes whenever $n<N$. Next, the fact that
\[(\alpha+\beta+2\theta+2\vartheta+N+k+3)_{N-k}
(\alpha+\beta+2\theta+2\vartheta+2n+4)_{k-n}\]
\[=\frac{(\alpha+\beta+2\theta+2\vartheta+2n+4)_{2N-2n-1}}
{(\alpha+\beta+2\theta+2\vartheta+k+n+4)_{N-n-1}}\] 
implies that we need only study the sum
\[\sum_{k=n}^{N}\frac{(-1)^{N-k}}{(N-k)!(k-n)!}\,
(\alpha+\beta+2\theta+2\vartheta+k+n+4)_{N-n-1}.\]
Shifting the sum by setting $j=k-n$ and $M=N-n$, introducing the variable 
\[\lambda=\alpha+\beta+2\theta+2\vartheta+2n+4,\]
and performing some manipulations, we find that the above sum transforms to
\[\frac{(-1)^M}{M!}\sum_{j=0}^M(-1)^j
{M\choose j}(\lambda+j)_{M-1}.\]
This is an iterated difference of order $M$, and as $(\lambda)_{M-1}$ is a 
polynomial of degree $M-1$, the iterated difference vanishes whenever 
$M>0$. The proof is complete.

\end{proof}

We now combine our results and obtain the norm expansion for the unit bidisk.

\begin{thm}
For any $f \in \Bergdiskgen$, we have 
\begin{equation*}
\|f\|^2_{\alpha,\beta,\theta,\vartheta}
=\sum_{N=0}^{+\infty}\frac{1}{\sigma(\alpha,\beta,\theta+N,\vartheta)}
\left|\left|\sum_{k=0}^{N}a_{k,N}\partial_{z_1}^{N-k}\oslash
\left[\partial_{z_1}^k f\right]\right|\right|^2_{\alpha+\beta+2\theta
+2\vartheta+2N+2},
%\label{normdisk}
\end{equation*}
where
\begin{multline*}
\frac{1}{\sigma(\alpha,\beta,\theta,\vartheta)}
=\frac{(\beta+1)\Gamma(\alpha+2)\Gamma(\theta+1)}
{(\alpha+\beta+2\theta+2\vartheta+3)\Gamma(\alpha+\theta+2)}\\
\times
{}_3F_2\left(\begin{array}{ccc}
\theta+1 & \alpha+\theta+\vartheta+2 & \alpha+\theta+\vartheta+2 
\\ {}&\alpha+\theta+2 & \alpha+\beta+2\theta+2\vartheta+4 
\end{array}\bigg|1\right) 
\end{multline*}
and
\begin{equation*}
a_{k,N}=\frac{(-1)^{N-k}}{k!(N-k)!}\frac{(\alpha+\theta+\vartheta+k+2)_{N-k}}
{(\alpha+\beta+2\theta+2\vartheta+N+k+3)_{N-k}}.
\end{equation*}
\end{thm}

\medskip

\noindent\bf The expression for the reproducing kernel of $\Bergdiskgen$. \rm
We now combine (\ref{kerrel2}), (\ref{kerneldecomp}) and the 
integral expression for $\Diffkernelz$ given in Lemma \ref{lemma-2}
and supply an explicit series and integral expression for the full 
reproducing kernel function of $\Bergdiskgen$. 

\begin{thm}
The reproducing kernel function of the space $\Bergdiskgen$ is
\begin{multline*}
\Kerneldiskgen(z,w)
=\sum_{N=0}^{+\infty}\sigma(\alpha,\beta,\theta+N,\vartheta)
(z_1-z_2)^N(\bar{w_1}-\bar{w_2})^N\\
\times\int_{\D}
\frac{\diff A_{\alpha+\beta+2\theta+2\vartheta+2N+2}(\xi)}
{[(1-\bar{\xi}z_1)(1-\bar{w_1}\xi)]^{\alpha+\theta+\vartheta+N+2}
[(1-\bar{\xi}z_2)(1-\bar{w_2}\xi)]^{\beta+\theta+\vartheta+N+2}}.
\end{multline*}
\label{thm-2}
\end{thm}
\medskip

\noindent\bf The weighted Hardy space case. \rm
We look at a special case of the identity of Theorem \ref{thm-2}. 
First, we set $\vartheta=0$ and note that in this case, the expression for 
the constant $\sigma(\alpha,\beta,\theta,\vartheta)$ reduces to
\[\frac{1}{\sigma(\alpha,\beta,\theta,0)}=
\frac{\Gamma(\alpha+2)\Gamma(\beta+2)\Gamma(\theta+1)\Gamma(\alpha+\beta+
2\theta+3)}{\Gamma(\alpha+\theta+2)
\Gamma(\beta+\theta+2)\Gamma(\alpha+\beta+\theta+3)},\]
and if we also put $\alpha=\beta=-1$, we get
\[\frac{1}{\sigma(-1,-1,\theta,0)}
=\frac{\Gamma(2\theta+2)}{(2\theta+1)[\Gamma(\theta+1)]^2}.\]

Next, we recall that in the limit $\alpha \rightarrow -1$, the weighted 
measure $\diff A_{\alpha}(z_1)$ degenerates to arc-length measure on the unit 
circle. This means that letting $\alpha$ and 
$\beta$ tend to $-1$ corresponds to considering the weighted Hardy space
$H^2_{\theta}(\D^2)$, with norm defined by
\begin{equation}
\|f\|^2_{H^2_{\theta}(\D^2)}=
\int_{\mathbb{T}^2}|f(z)|^2|z_1-z_2|^{2\theta}dm(z_1)dm(z_2),
\label{Hardybidisk}
\end{equation}
where $dm$ is the normalized Lebesgue measure on the unit circle.
Hence, Theorem \ref{thm-2} leads to a norm expansion for the weighted Hardy 
space. We state this as a corollary.

\begin{cor}
For any $f \in H^2_{\theta}(\D^2)$, we have
\begin{equation*}
\|f\|^2_{H^2_{\theta}(\D^2)}=\sum_{N=0}^{+\infty}
\frac{\Gamma(2\theta+2N+2)}{(2\theta+2N+1)[\Gamma(\theta+N+1)]^2}
\left|\left|\sum_{k=0}^N b_{k,N}\partial_{z_1}^{N-k}\oslash
\left[\partial_{z_1}^k f \right] \right|\right|^2_{2\theta+2N}
%\label{hardydiskexp},
\end{equation*}
where 
\begin{equation*}
b_{k,N}=\frac{(-1)^{N-k}}{k!(N-k)!}\frac{(\theta+k+1)_{N-k}}
{(2\theta+N+k+1)_{N-k}}.
\end{equation*}
\label{cor-hardydiskexp}
\end{cor}
\medskip

\section{Weighted Bergman spaces in the unit ball}
\label{Ball}

\noindent\bf Preliminaries. \rm
The unit ball in $\mathbb{C}^2$ is the set
\[\mathbb{B}^2=\big\{z=(z_1,z_2)\in\C^2:|z_1|^2+|z_2|^2<1\big\}.\]
For a survey of the function theory of the ball, we refer to \cite{Rud2};
see also \cite{Krantz} and \cite{Berg}.
We consider weighted spaces $\Lball$ consisting of (equivalence classes of)
Borel measurable functions $f$ on $\mathbb{B}^2$ with
\begin{equation*}
\|f\|^2_{\alpha, \beta,\theta}=
\int_{\mathbb B^2}|f(z)|^2\,|z_2|^{2\theta}
(1-|z_1|^2-|z_2|^2 )^\alpha
(1-|z_1|^2)^\beta\,\diff A(z_1,z_2) <+\infty,
%\label{lebspaceball}
\end{equation*} 
where $\alpha,\beta,\theta$ have 
$-1<\alpha,\beta,\theta<+\infty$,
and
$$\diff A(z_1,z_2)=\diff A(z_1)\diff A(z_2).$$
The Bergman spaces $\Aball$ is the  subspace of $\Lball$ consisting 
of functions $f$ that are holomorphic in $\mathbb B^2$. In this section, 
we find the orthogonal decomposition of functions in $\Aball$ along the
zero variety 
$$\big\{(z_1,z_2)\in \mathbb{B}^2:z_2=0\big\}.$$
As it turns out, the situation here is much easier to handle than in the 
bidisk case. 

By Taylor's formula, any function $f\in \Aball$ has a decomposition 
\[f(z)=\sum_{N=0}^{+\infty} g_N(z_1)\,z_2^N, \hspace{0.5cm}
\textrm{where} \hspace{0.5cm} g_N(z_1)=\frac 1{N!}\,
\partial_{z_2}^N f(z_1,0).\] 
It is easy to see that the summands in this decomposition are orthogonal in 
the space $\Aball$ for different $N$, and hence  
\begin{equation}
\|f\|^2_{\alpha,\beta,\theta}=\sum_{N=0}^{+\infty} 
\big\|g_N(z_1)z_2^N\big\|^2_{\alpha,\beta,\theta}. 
\label{derexpansion}
\end{equation}
\medskip

\noindent\bf The norm expansion for the ball. \rm
All we need is the following lemma.

\begin{lem}
We have that
\begin{equation*}
\big\|g_N(z_1)z_2^N\big\|^2_{\alpha,\beta,\theta}=
\frac {\Gamma(\alpha+1)\Gamma(\theta+N+1)}
{(\alpha+\beta+\theta+N+2)\Gamma(\alpha+\theta+N+2)}
\,\|g_N\|_{\alpha+\beta+\theta+N+1}^2. 
\end{equation*}
\label{lemma-ballemb}\end{lem}

\begin{proof}
We make the change of variables 
\[z_1=z_1,\quad z_2=(1-|z_1|)^{1/2}u,\]
and get
\begin{multline}
\big\|g_N(z_1)z_2^N\big\|^2_{\alpha,\beta,\theta}=
\int_{\mathbb B^2}|g_N(z_1)|^2\,|z_2|^{2\theta+2N}
(1-|z_1|^2-|z_2|^2)^\alpha
(1-|z_1|^2)^\beta\, \diff A(z_1,z_2)\\=
\int_{\mathbb D}|g_N(z_1)|^2 (1-|z_1|^2)^{\alpha+\beta+\theta+N+1} 
\diff A(z_1) \times
\int_{\mathbb D}|u|^{2(\theta+N)} (1-|u|^2)^\alpha \diff A(u),
\label{ballembpr}
\end{multline}
whence the assertion follows.
\end{proof}

Finally, we obtain 

\begin{thm}
\label{ball-t1}
For any $f\in\Aball$, 
\begin{equation*}
\|f\|^2_{\alpha,\beta,\theta}=
\sum_{N=0}^{+\infty} 
\frac{\Gamma(\alpha+1)\Gamma(\theta+N+1)}
{(\alpha+\beta+\theta+N+2)\Gamma(\alpha+\theta+N+2)(N!)^2}  
\,\big\|\partial_{z_2}^N f(z_1,0)  
\big\|^2_{\alpha+\beta+\theta+N+1}. 
\label{normball}
\end{equation*}
\end{thm}
\medskip

\noindent\bf Weighted Hardy spaces. \rm
As in the case of the bidisk, we derive a corollary concerning 
weighted Hardy spaces also for the ball. 
We have
$$\lim_{\alpha\to -1+0}\, (\alpha+1)(\alpha+2)
\|f\|_{\alpha,\beta,\theta}^2 = 
\int_{\partial \mathbb B^2}|f(z)|^2\,|z_2|^{2\theta}
(1-|z_1|^2)^\beta\,
d\sigma(z),$$ 
where $d\sigma$ is the normalized Lebesgue measure on $\partial \mathbb B^2$. 
The right hand side of the last formula represents the norm of a function 
in  the weighted Hardy space denoted by $H^2_{\beta,\theta}(\mathbb B^2)$. 
As a corollary, we obtain the following decomposition of the norm of 
functions from this weighted Hardy space: 
\begin{cor}
\label{ball-cor1}
For any $f\in H^2_{\beta,\theta}(\mathbb B^2)$, 
$$\|f\|^2_{H^2_{\beta,\theta}(\mathbb B^2)} = 
\sum_{N=0}^{+\infty} \frac 1{(\beta+\theta+N +1)(N!)^2}\,
\big\|\partial_{z_2}^N f(z_1,0)  
\big\|^2_{\beta+\theta+N}. $$
\end{cor}
\medskip

\noindent \bf An expression for the reproducing kernel of $\Aball$. \rm
Now, we derive an explicit formula for the reproducing kernel for the 
space $\Aball$. In conformity with the notation in the introduction, 
we denote by $\Jball$ the subspace of $\Aball$ consisting of functions
of the form $f(z)=z_2^Ng(z_1)$. An easy calculation based on Lemma 
\ref{lemma-ballemb} establishes the following result.

\begin{lem}
The reproducing kernel  
for $\Jball$ is given by the formula
\begin{equation*}
Q_{\alpha,\beta,\theta,N}\big(z,w\big) \\
=\frac{(\alpha+\beta+\theta+N+2)\Gamma(\alpha+\theta+N+2)}
{\Gamma(\alpha+1)\Gamma(\theta+N+1)}\times
\frac{(z_2\bar w_2)^N}{(1-z_1\bar w_1)^{\alpha+\beta+\theta+N+3}}. 
\end{equation*}
\label{lemma-Qkernelball}
\end{lem}

Since $\Aball$ is the orthogonal sum of the subspaces 
$\Jball$, its reproducing kernel $P_{\alpha,\beta,\theta}$ is given by 
the sum 
\begin{multline*}
P_{\alpha,\beta,\theta}(z,w)=
\sum_{N=0}^{+\infty}
Q_{\alpha,\beta,\theta,N}(z,w)\\
=\frac{\Gamma(\alpha+\theta+2)}{\Gamma(\alpha+1)\Gamma(\theta+1)}
 \frac {1}{(1-z_1\bar w_1)^{\alpha+\beta+\theta+3}} 
\sum_{N=0}^{+\infty} \frac{(\alpha+\beta+\theta+N+2)(\alpha+\theta+2)_N}
{(\theta+1)_N}\left(\frac{z_2\bar w_2}{1-z_1\bar w_1}\right)^N\\ = 
\frac{\Gamma(\alpha+\theta+2)}{\Gamma(\alpha+1)\Gamma(\theta+1)} 
\frac 1{(1-z_1\bar w_1)^{\theta+\alpha+\beta+3}}\\ \times
\left[(\alpha+\theta+2)\,\,{}_2F_1\left(\alpha+\theta+3,1;\theta+1;
\frac{z_2\bar w_2}{1-z_1\bar w_1}\right)\right. + \left.
\beta\,\, {}_2F_1\left(\alpha+\theta+2,1;\theta+1;
\frac{z_2\bar w_2}{1-z_1\bar w_1}\right) \right]. 
\end{multline*}
Here, ${}_2F_1$ stands for the classical Gauss hypergeometric function. 
We formulate the result as a theorem.

\begin{thm}
The kernel function for the space $\Aball$ is
\begin{multline}
P_{\alpha,\beta,\theta}(z,w)
=\frac{\Gamma(\alpha+\theta+2)}{\Gamma(\alpha+1)\Gamma(\theta+1)}
\frac{1}{(1-\bar{w_1}z_1)^{\alpha+\beta+\theta+3}}\\
\times
\left[(\alpha+\theta+2)\,\,{}_2F_1\left(\alpha+\theta+3,1;\theta+1;
\frac{z_2\bar{w_2}}{1-\bar{w_1}z_1}\right) \right.\left.+
\beta\,\, {}_2F_1\left(\alpha+\theta+2,1;\theta+1;\frac{z_2\bar{w_2}}
{1-\bar{w_1}z_1}\right)\right].
\end{multline}
\end{thm}

\begin{rem} It would be natural to consider more general Hilbert space norms
of the type
\begin{equation*}
\|f\|^2_{\alpha, \beta,\theta,\gamma}=
\int_{\mathbb B^2}|f(z)|^2\,|z_2|^{2\theta}
(1-|z_1|^2-|z_2|^2 )^\alpha
(1-|z_1|^2)^\beta\,(1-|z_2|^2)^\gamma\,\diff A(z_1,z_2),
\end{equation*} 
which are symmetric with respect to an interchange of the variables $z_1$ and
$z_2$ (if simultaneously $\beta$ and $\gamma$ are interchanged).  
Here, we must suppose that $-1<\alpha,\beta,\gamma,\theta<+\infty$. 
The already treated case corresponds to $\gamma=0$. 
The above analysis applies here as well, but, unfortunately, the formulas 
become rather complicated; this is why we work things out for $\gamma=0$ only.
\end{rem}

\section{Weighted Bargmann-Fock spaces in $\C^2$}
\label{Bargmann}

\noindent\bf Preliminaries. \rm
Fix a real parameter $\gamma$ with $0<\gamma<+\infty$. The classical 
one-variable Bargmann-Fock space -- denoted by 
$A^2_\gamma(\C)$ -- consists of all entire functions of one complex variable 
with
\begin{equation}
\|f\|^2_{\gamma}=
\int_\C|f(z)|^2\,e^{-\gamma|z|^2}\diff A(z)<+\infty,
\label{eq-A2gamma}
\end{equation}
the associated sesquilinear inner product is denoted by 
$\langle\cdot,\cdot\rangle_\gamma$. 
The reproducing kernel of this Hilbert space is well-known:
$$P_\gamma(z,w)=\gamma\,e^{\gamma\bar w z},\qquad z,w\in\C.$$

Next, fix real parameters $\alpha,\beta,\theta$ with 
$0<\alpha,\beta<+\infty$ and $-1<\theta<+\infty$.
We consider the Hilbert space 
$\LBargmann$ of
all (equivalence classes of) Borel measurable functions $f$ in the bidisk 
subject to the norm boundedness condition
\begin{equation*}
\|f\|^2_{\alpha,\beta,\theta}=
\int_\C\int_{\C}|f(z_1,z_2)|^2|z_1-z_2|^{2\theta}\,
e^{-\alpha|z_1|^2-\beta|z_2|^2}\diff A(z_1)\diff A(z_2)<+\infty;
%\label{lebspaceB}
\end{equation*}
we let $\langle\cdot,\cdot\rangle_{\alpha,\beta,\theta}$ denote
the associated sesquilinear inner product. 
The {\em weighted Bargmann-Fock space} $\Bargmann$ is the subspace of 
$\LBargmann$ consisting of the entire functions.

The reproducing kernel for $\Bargmann$ will be denoted by
\[\KernelB=
\KernelB(z,w),\]
where $z=(z_1,z_2)$ and $w=(w_1,w_2)$ are two points in $\C^2$. The kernel 
defines an orthogonal projection of the space $\LBargmann$ onto the weighted 
Bargmann-Fock space via the formula
\begin{multline*}
\KernelB[f](z)=\big\langle f,
\KernelB
(\cdot,z)\big\rangle_{\alpha,\beta,\theta}\\
=\int_{\C}\int_{\C} f(w)P_{\alpha,\beta,\theta}
(z,w)\,|w_1-w_2|^{2\theta}\,e^{-\alpha|z_1|^2-\beta|z_2|^2}
\diff A(w_1)\diff A(w_2);
%\label{reprodpropB}
\end{multline*}
as indicated, we shall write $\KernelB[f]$ for the projection of a 
function $f\in\LBargmann$.

In the case $\theta=0$, the reproducing kernel it is readily computed:
\begin{equation*}
P_{\alpha,\beta,0,0}(z,w)=
\alpha\beta\,e^{\alpha z_1\bar w_1+\beta z_2\bar w_2}.
\end{equation*}

We consider the polynomial $p(z_1,z_2)=z_1-z_2$ in the context of the 
introduction. In particular, for non-negative integers $N$, we consider 
the subspaces $\ZeroNB$ of functions in $\Bargmann$ that vanish up to 
order $N$ along the diagonal 
$$\diag(\C)=\big\{(z_1,z_2)\in \C^2:z_1=z_2\big\}.$$
Being closed subspaces of a reproducing kernel space, the spaces $\ZeroNB$ 
possess reproducing kernels of their own. We shall write
\[\KernelBNB=\KernelBNB(z,w)\]
for these kernel functions.
The operators associated with the kernels project the space
$L^2_{\alpha,\beta,\theta}(\C^2)$
orthogonally onto $\ZeroN$. As before, we write
$\KernelBNB[f]$ for the projection of a function.

Next, we define the spaces $\DiffNB$ by setting
\[\DiffNB=\ZeroNB \ominus \ZeroNpB.\]
The spaces $\DiffNB$ also admit reproducing kernels, and their
kernel functions are of the form
\begin{equation*}\DiffkernelN(z,w)
=\KernelBNB(z,w)
-P_{\alpha,\beta,\theta,N+1}(z,w).
\end{equation*}
We shall write $\DiffkernelzB$ for the kernel 
$Q_{\alpha,\beta,\theta,0}$.
\medskip

As in the case of the weighted Bergman spaces on the bidisk, we make the 
following observation. We suppress the proof, as it is virtually identical to
that of Lemma \ref{lemma-1}.

\begin{lem}
We have
\begin{equation*}
\KernelBNB(z,w)
=(z_1-z_2)^N(\bar{w_1}-\bar{w_2})^N
\KernelNB(z,w)
%\label{kerrel1B}
\end{equation*}
for $z,w\in\C^2$. 
\label{lemma-10}
\end{lem}

If we write, as in the introduction, 
$$\hispace(\C^2)=A^2_{\alpha,\beta,\theta}(\C^2),$$
we may identify the spaces $\hispace_N(\C^2)$,
$$\hispace_N(\C^2)=A^2_{\alpha,\beta,\theta+N}(\C^2),
\qquad N=0,1,2,\ldots,$$
and the spaces $\gispace_N(\C^2)$ as well:
$$\gispace_N(\C^2)=
\mospace_{\alpha,\beta,\theta+N,0}(\C^2),\qquad N=0,1,2,\ldots.$$
As a consequence, we get that
\begin{equation}
\DiffkernelNB(z,w)
=(z_1-z_2)^N(\bar{w_1}-\bar{w_2})^N
\DiffkernelNpNB(z,w).
\label{kerrel2B}
\end{equation}
By (\ref{exp-1}), we have the kernel function expansion
\begin{equation}
\KernelB(z,w)
=\sum_{N=0}^{+\infty}\DiffkernelNB(z,w),
\label{kerneldecompB}
\end{equation}
while the orthogonal norm expansion (\ref{norm-1}) reads
\begin{equation}
\|f\|^2_{\alpha,\beta,\theta}
=\sum_{N=0}^{+\infty}\big\|\DiffkernelNB[f]
\big\|^2_{\alpha,\beta,\theta},
\qquad f\in\Bargmann. 
\label{normdecompB}
\end{equation}

Our next objective is to identify the Hilbert space of restrictions to the
diagonal of 
$\hispace_N(\C^2)=A^2_{\alpha,\beta,\theta+N}(\C^2)$, as
well as to calculate the reproducing kernel of $\hispace_N(\C^2)$ on the
set $\C^2\times\diag(\C)$.
\medskip

\noindent\bf Unitary operators.
\rm 
The rotation operator $R_{\phi}$ (for a real parameter $\phi$) defined for 
$f\in\Bargmann$ by
\[R_{\phi}[f](z_1,z_2)=f(e^{i\phi}z_1,e^{i\phi}z_2)\]
is clearly unitary, and we shall make use of it shortly. 
The following lemma supplies us with yet another family of unitary operators.

\begin{prop}
For every $\lambda\in\C$, the operator
\begin{equation*}
U_{\lambda}[f](z_1,z_2)=
e^{-(\alpha+\beta)|\lambda|^2/2}
e^{-\alpha\bar\lambda z_1-\beta\bar\lambda z_2}
\,f\left(z_1+\lambda,z_2+\lambda\right)
%\label{unitopB}
\end{equation*}
is unitary on the space $\Bargmann$, and its adjoint is $U_\lambda^*=
U_{-\lambda}$. 
\end{prop}

The proof amounts to making a couple of elementary changes of variables in 
integrals, and is therefore left out.
\medskip

\noindent\bf The reproducing kernel on the diagonal. \rm
We now use the operators $R_{\phi}$ and $U_{\lambda}$ to 
compute the reproducing kernel on the set $\C^2\times\diag(\C)$.

\begin{thm}
We have
\begin{equation*}
\KernelB((z_1,z_2),(w_1,w_1))=\DiffkernelzB((z_1,z_2),(w_1,w_1))=
\sigma(\alpha,\beta,\theta)\,e^{\alpha\bar w_1 z_1+\beta\bar w_1 z_2}.
%\label{diagker}
\end{equation*}
for $z_1,z_2,w_1\in\C$. Here, $\sigma(\alpha,\beta,\theta)$ is the 
positive constant  given by 
\begin{equation*}
\frac{1}{\sigma(\alpha,\beta,\theta)}
=\int_\C\int_{\C}|z_1-z_2|^{2\theta}e^{-\alpha|z_1|^2-\beta|z_2|^2}
\diff A(z_1)\diff A(z_2)=
\frac{(\alpha+\beta)^\theta}{(\alpha\beta)^{\theta+1}}\,\Gamma(\theta+1).
%\label{sigma}
\end{equation*}
\label{thm-10}
\end{thm}

\begin{proof}
By using the unitarity of the rotation operator $R_{\phi}$, we get as in
the proof of Theorem \ref{thm-1} that the function $\KernelB(\cdot,0)$
is positive constant, which we denote by $\sigma(\alpha,\beta,\theta)$. 

Next, take $\lambda \in \C$ and $f \in \Bargmann$.
As the operators $U_{\lambda}$ are unitary, and as 
$U_{\lambda}^*=U_{-\lambda}$, we find that
\begin{multline*}
e^{-(\alpha+\beta)|\lambda|^2/2}\,
f(\lambda,\lambda)
=U_{\lambda}[f](0)
=\big\langle U_{\lambda}[f],
\KernelB(\cdot,0)
\big\rangle_{\alpha,\beta,\theta}\\
=\big\langle f, U_{-\lambda}[\KernelB
(\cdot,0)]\big\rangle_{\alpha,\beta,\theta}
=\sigma(\alpha,\beta,\theta)\,
\big\langle f,U_{-\lambda}[1]\big\rangle_{\alpha,\beta,\theta}.
\end{multline*}
This equality together with the uniqueness of reproducing kernels
establishes that
\[\KernelB((z_1,z_2),(\lambda,\lambda))
=\sigma(\alpha,\beta,\theta)\,e^{(\alpha+\beta)|\lambda|^2/2}\,
U_{-\lambda}[1](z_1,z_2),\]
which is the desired result.
The explicit expression for the constant in terms of an integral over the 
bidisk follows if we apply the reproducing property of the kernel applied 
to the constant function $1$. The evaluation of the integral in terms
of the Gamma function is done by performing a suitable change of
variables.
\end{proof}
\medskip

\noindent\bf Restrictions of reproducing kernels. \rm
From the previous subsection, we have that 
\[\KernelB((z_1,z_2),(w_1,w_1))
=\sigma(\alpha,\beta,\theta)
\,e^{\alpha\bar w_1 z_1+\beta\bar w_1 z_2}.\]
For continuous functions 
$f \in L^2_{\alpha,\beta,\theta}(\C^2)$, 
we use the notation $\oslash f$ for the restriction 
to the diagonal of the function, that is,
\[(\oslash f) (z_1)=f(z_1,z_1),\qquad z_1\in\C,\]
just like in Section \ref{Bidisk}. We fix $w_1$ and 
apply this operation to the reproducing kernel function of $\Bargmann$. 
We obtain
\[\oslash\, \KernelB(z_1,(w_1,w_1))
=\sigma(\alpha,\beta,\theta)\,e^{(\alpha+\beta)\bar w_1 z_1}.\]
and we see that the restriction of the kernel
coincides with a multiple of the reproducing kernel function for
the space $A^2_{\alpha+\beta}(\C)$.
By the theory of reproducing kernels (see \cite{Sa}), this means that 
the induced norm for the space 
$\mospace_{\alpha,\beta,\theta,0}(\C^2)$ coincides with a 
multiple of the norm in the aforementioned Bargmann-Fock space of one 
variable. An immediate consequence of this fact is the inequality
\begin{equation}
\frac{\alpha+\beta}{\sigma(\alpha,\beta,\theta)}
\|\oslash f\|^2_{\alpha+\beta}
\leq \|f\|^2_{\alpha,\beta,\theta},\qquad f\in\Bargmann,
\label{rnorineqB}
\end{equation}
and, more importantly, the equality
\begin{equation}
\frac{\alpha+\beta}{\sigma(\alpha,\beta,\theta)}\,
\|\oslash f\|^2_{\alpha+\beta}
=\big\|\DiffkernelzB[f]\big\|^2_{\alpha,\beta,\theta},\qquad
f\in\mospace_{\alpha,\beta,\theta,0}(\C^2).
\label{rnormeqB}
\end{equation}
The notation on the left hand sides of (\ref{rnorineqB}) and (\ref{rnormeqB})
is in conformity with (\ref{eq-A2gamma}).

The next step in our program is to compute the kernel function 
$\DiffkernelzB$. 

\begin{prop}
The kernel function for the space 
$\mospace_{\alpha,\beta,\theta,0}(\C^2)$
is given by 
\begin{equation*}
\DiffkernelzB(z,w)
=\frac{(\alpha\beta)^{\theta+1}}{(\alpha+\beta)^\theta\Gamma(\theta+1)}
\,e^{(\alpha\bar w_1+\beta\bar w_2)(\alpha z_1+\beta z_2)/(\alpha+\beta)},
\qquad z,w\in \C^2.
\end{equation*}
\label{prop-11}
\end{prop}

\begin{proof}
In the notation of the introduction, we have the identity (with $N=0$)
$$k^{\gispace_{N}(\Omega)}(z,w)=
\big\langle \oslash_p k^{\hispace_{N}(\Omega)}_{w},
\oslash_p k^{\hispace_{N}(\Omega)}_{z}
\big\rangle_{\hispace_{N}(\zeroset_p)},$$
by a combination of (\ref{def-l}) and (\ref{eq-l}). In the notation of this
section, this means that
\begin{equation*}
\DiffkernelzB(z,w)
=\frac{\alpha+\beta}{\sigma(\alpha,\beta,\theta)}
\,\big\langle \oslash\KernelB(\cdot,w),\oslash\KernelB(\cdot,z)\big\rangle
_{\alpha+\beta},
\qquad z,w\in \C^2,
\end{equation*}
so that by applying Theorem \ref{thm-10}, we get 
\begin{equation*}
\DiffkernelzB(z,w)
=(\alpha+\beta)\sigma(\alpha,\beta,\theta)\int_\C 
e^{(\alpha z_1+\beta z_2)\bar\xi}\,e^{(\alpha\bar w_1+\beta\bar w_2)\xi}\,
e^{-(\alpha+\beta)|\xi|^2}\,\diff A(\xi).
\end{equation*}
It just remains to evaluate the integral.
\end{proof}
\medskip

\noindent\bf The expression for the reproducing kernel of $\Bargmann$. \rm
In view of Proposition \ref{prop-11}, (\ref{kerrel2B}), and 
(\ref{kerneldecompB}), we may now derive an explicit expression for the 
reproducing kernel of $\Bargmann$.

\begin{cor}
The reproducing kernel for $\Bargmann$ is given by
$$\KernelB(z,w)=\frac{(\alpha\beta)^{\theta+1}}{(\alpha+\beta)^\theta}
\,e^{(\alpha\bar w_1+\beta\bar w_2)(\alpha z_1+\beta z_2)/(\alpha+\beta)}
\,E_\theta\left(\frac{\alpha\beta}{\alpha+\beta}\,
(z_1-z_2)(\bar w_1-\bar w_2)\right),$$
where
$$E_\theta(x)=\sum_{N=0}^{+\infty}\frac{x^N}{\Gamma(\theta+N+1)},\qquad
x\in\C.$$
\end{cor}
\medskip

\noindent\bf The diagonal norm expansion for the Bargmann-Fock space. \rm
Having obtained the reproducing kernel in explicit form, we only need to
write the norm decomposition (\ref{normdecompB}) in desired form.

\begin{lem}
For each $N=0,1,2,\ldots$, we have the equality of norms 
\begin{equation*}
\big\|\DiffkernelNB[f]\big\|^2_{\alpha,\beta,\theta}
=\frac{(\alpha+\beta)^{\theta+N+1}\Gamma(\theta+N+1)}
{(\alpha\beta)^{\theta+N+1}}
\left|\left|\oslash \left[\frac{P_{\alpha,\beta,\theta,N}[f]}
{(z_1-z_2)^N}\right]\right|\right|^2_{\alpha+\beta},
\end{equation*}
for all $f\in\Bargmann$.
\label{lemma-12}
\end{lem}

\begin{proof}
The statement follows from a combination of Lemma \ref{lemma-10} and 
(\ref{rnormeqB}), plus the evaluation of $\sigma(\alpha,\beta,\theta+N)$.
\end{proof}

All that remains for us to do is to make the right hand side of the
expression in Lemma \ref{lemma-12} sufficiently explicit.

\begin{lem}
For all $k=0,1,2,\ldots$, and each $f\in\Bargmann$, we have
\begin{equation}
\oslash\, \partial_{z_1}^{k}[f]=
\sum_{n=0}^{k}n!{k\choose n}
\left(\frac{\alpha}{\alpha+\beta}\right)^{k-n}\,
\partial_{z_1}^{k-n}
\oslash \left[\frac{P_{\alpha,\beta,\theta,n}[f]}{(z_1-z_2)^n}
\right].
\end{equation}
\label{lemma-13}
\end{lem}

\begin{proof}
We observe that 
$$\oslash\,[\partial_{z_1}^j\DiffkernelNB](z_1,(w_1,w_2))=
\left(\frac{\alpha}{\alpha+\beta}\right)^j\,
\partial_{z_1}^j\oslash\,[\DiffkernelNB](z_1,(w_1,w_2)).$$
The rest of the proof is obtained by mimicking the arguments of Lemma 
\ref{lemma-3}.
\end{proof}
\medskip

It is quite easy to invert Lemma \ref{lemma-13}:

\begin{lem}
for all $N=0,1,2,\ldots$ and each $f\in\Bargmann$, the equality
\begin{equation}
\oslash\left[ \frac{P_{\alpha,\beta,\theta}[f]}{(z_1-z_2)^N}
\right](z_1) 
=\sum_{k=0}^{N}\frac{(-1)^{N-k}}{k!(N-k)!}\,
\left(\frac{\alpha}{\alpha+\beta}\right)^{N-k}\,
\partial_{z_1}^{N-k}\oslash \left[
\partial_{z_1}^k f\right](z_1),\qquad z_1\in\C,
\end{equation}
holds for each $f \in \Bargmann$.
\end{lem}

\begin{proof}
In view of the Lemma \ref{lemma-13}, it is enough to check that
\[\sum_{k=n}^N\frac{(-1)^{N-k}}{k!(N-k)!}\,\,
n!{k\choose n}
\left(\frac{\alpha}{\alpha+\beta}\right)^{N-n}=\delta_{n,N},\]
where $\delta_{n,N}$ is the Kronecker delta. Firstly,
we observe that the equality holds for $n=N$. Secondly, 
we observe that it is equivalent to show that
$$\sum_{k=n}^N\frac{(-1)^{N-k}}{k!(N-k)!}\,\,
n!{k\choose n}=\sum_{k=n}^N\frac{(-1)^{N-k}}{(N-k)!(k-n)!}=0$$
whenever $n<N$. The expression in the middle is an $(N-n)$-th 
difference of a constant function, which of course is $0$ for $n<N$.
We are done. 
\end{proof}

We now combine our results and obtain the norm expansion for the 
Bargmann-Fock space.

\begin{thm}
Let $c_{k,N}(\alpha,\beta)$ be given by
$$c_{k,N}(\alpha,\beta)=(-1)^{N-k}{N\choose k}
\left(\frac{\alpha}{\alpha+\beta}\right)^{N-k}.$$
Then, for each $f\in\Bargmann$, we have 
\begin{equation*}
\|f\|^2_{\alpha,\beta,\theta}
=\sum_{N=0}^{+\infty}\frac{(\alpha+\beta)^{\theta+N+1}\Gamma(\theta+N+1)}
{(\alpha\beta)^{\theta+N+1}[N!]^2}
\left|\left|\sum_{k=0}^{N}c_{k,N}(\alpha,\beta)\,\partial_{z_1}^{N-k}\oslash
\left[\partial_{z_1}^k f\right]\right|\right|^2_{\alpha+\beta}.
\label{normdisk}
\end{equation*}
\end{thm}

\begin{rem}
There is an alternative way to obtain the norm expansion and the explicit 
expression for the reproducing kernel in the Bargmann-Fock space $\Bargmann$. 
The change of variables 
$$\left\{
\begin{array}{lll}z_1&=&w_1+\beta w_2 \\ z_2&=&w_2-\alpha w_2 \end{array}
\right. 
$$
transforms the norm in  $\Bargmann$ into the expression
$$\|f\|_{\alpha,\beta,\theta}=(\alpha+\beta)^{2\theta+2} \int_\C\int_\C
|g(w_1,w_2)|^2\,|w_2|^{2\theta}
e^{-(\alpha+\beta)|w_1|^2-\alpha\beta(\alpha+\beta)|w_2|^2}\,
dA(w_1)dA(w_2),
$$
where $g(w_1,w_2)=f(w_1+\beta w_2,w_1-\alpha w_2)$. 
The reproducing kernel and the norm expansion about hyperplane $w_2=0$ with 
respect to the latter norm can be calculated by separation of variables.
Shifting back to the original variables $(z_1,z_2)$, then, we obtain the 
reproducing kernel and norm expansion for $\Bargmann$.  
\end{rem}

\end{document}